\documentclass[12pt]{article}
\usepackage{amsmath,amssymb,amsfonts,amsthm,graphicx}
\usepackage{enumerate}
\usepackage[mathscr]{eucal}
\usepackage{color} 
\title{A new stochastic order based on discrete Laplace transform and some ordering results of the order statistics  }
\author{ Fatemeh Gharari and Masoud Ganji
\hspace{.5cm} 
\\
\footnotesize{Department of Statistics and Computer Science, University of Mohaghegh Ardabili, Ardabil, IRAN} \\
\footnotesize{E-mails:  f.gharari@uma.ac.ir, mganji@uma.ac.ir  }\\
}
\date{}
\begin{document}
\maketitle \thispagestyle{empty}
\newtheorem{theorem}{Theorem}[section]
\newtheorem{lemma}[theorem]{Lemma}
\newtheorem{proposition}[theorem]{Proposition}
\newtheorem{corollary}[theorem]{Corollary}
\theoremstyle{definition}
\newtheorem{definition}[theorem]{Definition}
\newtheorem{example}[theorem]{Example}
\newtheorem{xca}[theorem]{Exercise}
\theoremstyle{remark}
\newtheorem{remark}[theorem]{Remark}
\numberwithin{equation}{section}
\begin{abstract}
This paper aims to study a new stochastic order based  upon discrete Laplace transforms. By this order, in a setup  where the 
sample size is  random, having   discrete  delta and  nabla   distributions,  we obtain some ordering results involving ordinary and fractional order statistics. Some applications in frailty models and   reliability  are presented as well.
\end{abstract}
\textbf{Keywords:}  Fractional derivative, Laplace transform, Order statistics, Time scale, Stochastic order.  \\

\section{Introduction}
Let $X_{1}, ..., X_{n}$ be independent and identically distributed (iid) random variables from a
distribution $F. $ The $ i $th order statistics (os) will be denoted by ${X}_{i:n},\,i=1,2,...,n. $ For further details about os, we refer the reader to Arnold et al. (1992).
\\
In recent years, many researchers have studied stochastic comparisons of os 
  for iid random variables  when the sample size is fixed.   Shaked
and Wong (1997) have studied stochastic comparisons of os, ${X}_{1:n}$ and
${X}_{n:n}.$  Some other results regarding stochastic comparisons of the os
 ${X}_{1:n}$  and ${X}_{n:n}$  can be obtained in Bartoszewicz (2001).
Nanda et al. (2005) and Nanda and Shaked  (2008)  generalized the results proved in
Shaked and Wong (1997) and Bartoszewicz (2001) to os other than the minimum and
the maximum.  They have studied the closure of the different partial orders under the
formation of a  $r$-out-of-$N$ system when the number of components $N,$ forming the
system, is a random variable with support ${k, k + 1, ... },$ for some fixed positive
integers $k. $ An open problem in stochastic order theory is \textit{stochastic comparisons of fractional order statistics} (fos) when the sample size is random.
\\
The concept of fos was introduced
by Stigler (1977) via a Dirichlet process definition. This type of os was developed as a technical
device which may be viewed as defining a continuum of os for any sample
size.  
There are two principle definitions for fos .The first definition presented by Jones (2002)   as a random convex combination of
${X}_{i:n}$ and ${X}_{i+1:n},$ namely
\begin{equation*}                                                                                                                                                                      {X}_{a:n}={X}_{i+c:n}=(1-c){X}_{i:n}+c{X}_{i+1:n},
\end{equation*}
where $1<a<n,\,\,i=\lfloor a\rfloor$ denotes the largest integer less than or equal to $ a,\, c=a-i,$ and ${X}_{i:n}$ is the $i $th os   from an independent and identically distributed sample size $ n$ from distribution $ F.$
The second definition of fos is due to Stigler (1977). Let ${X}_{a:n}= F^{-1}(U_{a:n}),$ where $U_{a:n}$ is  fractional associated with the uniform (on $[0,1]$)  defined to be a quantity following the beta distribution (on $[0,1]$) with parameters $ a $ and $ n-a+1, $ written $ Beta(a, n-a+1).$ This definition  is appropriate for $1<a<n+1.$ As an application of fos,  non-parametric confidence intervals for quantiles  can be calculated using
fos 
(see Hutson (1999)). Another application is to construct approximate prediction intervals of os based on fos (see Basiri et al.(2016)). As noted by Hutson (1999), in general, fos cannot
be calculated from the sample, since their indices  do not need to take integer values.  An open problem in order statistics theory is \textit{finding the distribution of a fos}. 
\\
In this study, we try to find solutions for these two problems. To this end, first we present   new definitions of os and fos by using discrete
fractional calculus. On the other hand, we know that the
existence of the  delta and nabla discrete distributions were proved by using discrete
fractional calculus (see Ganji and Gharari (2018)). We see that these two definition processes were generated on a construction equipment of delta and nabla
calculus. Considering this important point, we can find solutions by looking for some relationships
 between os and fos with  nabla and delta types of probability distributions. In this regard, some of important constructions will form the basis of our theory.
First, we suppose that the  sample size $ N $ is   a random variable having discrete nabla  or delta distribution. Under this situation, we give two general expressions for the  probability density functions (pdfs) of the $ r$th os and fos in terms of  the $r$th derivatives of the Laplace transform of  $ N.$ For the next problem, we define a new stochastic order based on the ordinary and fractional derivatives  of the Laplace transform of the sample size $ N $  and compare os and fos in a setup where the sample size is random.
\\
This study has two main advantages in comparison with previous works regarding  the comparison of os  when the sample size is random.
First,  the new stochastic order appears  the relation between the $i$th os  and the $i$th  derivative  of the Laplace transform of the  random sample size $ N,$ where $ i$ and $ N\in\mathbb{N}_{1}. $  Therefore, we present a new method for the comparison of   os other than the minimum and
the maximum. On the other hand, we  generalize   the new stochastic order  for the comparison of  fos when the sample size is random. Similarly in this case, the new stochastic order reveals the relationship between  $\gamma$th fos, $n-1<\gamma\leq n,\,n\in \mathbb{N}$  when the sample size $ N $ is random,  and the $\gamma$th fractional derivative  of the Laplace transform of $ N ,$ where $ N\in\mathbb{N}_{\alpha-1},\,\alpha-1>0.$
\\
The article  is organized as follows. The second section contains  some
definitions and preliminary results that are needed for the  main results to be established in the next sections.
  The third section is essential for starting our theory and contains a new definition of a fos by spacings. In  the fourth section, we give another definition of a  fos by fractional operators and there represent two general expressions  for the  pdfs   of the   $r$th os and fos when
the sample  size $N$  is random.  Also, we introduce two new stochastic orderings based on both  fractional and  ordinary $r$th derivations of the  Laplace transforms of $N.$ At the next section, some stochastic orders based on discrete Laplace transforms
  are defined along with some results and properties. Some of their applications are presented in the final section.


\section{Preliminaries}
For  a non-negative random variable $X$ with df  $F,$ its Laplace transform is defined as
\begin{equation*}
\mathcal{L}_{X}(s)=\int_{0}^{\infty}e^{-sx}dF(x).
\end{equation*}
 Let $X$ and $Y$ be two non-negative random variables. $X$ is
said to be smaller than $Y$ in the Laplace transform order (denoted by $X\leqslant_{lt}Y$) if $ L_{Y} (s)\leqslant L_{X}(s)$ for all $s>0.$ Shaked and Wong (1997)
defined two stochastic orders based on ratios of the Laplace transforms: $X$ is said to be smaller than $Y$ in the Laplace transform ratio
order (denoted by $X\leqslant_{lt-r}Y$) if $ L_{Y} (s)/L_{X}(s)$ is decreasing in $s>0;$ and $X$ is said to be smaller than $Y$ in the reversed Laplace transform
ratio order (denoted by $ X\leqslant_{ r-lt-r}Y$) if $[1- L_{Y} (s)]/[1- L_{X}(s)]$ is decreasing in $s>0.$   Li et al. (2009) defined a stochastic order based on ratio of differentiated Laplace transforms: $X$  is said to be smaller than $Y$ in the differentiated Laplace transform ratio order (denoted by $X \leqslant_{d-lt-r}Y$) if 
$L^{^{\prime}}_{Y} (s)/L^{^{\prime}}_{X}(s)$
is decreasing in $s>0.$
Now, let us introduce a new order which, as will be seen in Section 6, has some potential applications.
\begin{definition}
$X$  is said to be smaller than $Y$ in the $i-$differentiated Laplace transform ratio order (denoted by $X \leqslant_{d^{^{(i)}}-lt-r}Y$) if 
$L^{^{(i)}}_{Y} (s)/L^{^{(i)}}_{X}(s)$
is decreasing in $s>0.$
\end{definition}
In the following, we generalize these stochastic orders for discrete Laplace transforms, which, as will be seen in the next sections, have some potential applications.\\
For real numbers  $ a $ and  $ b, $ we denote $ \mathbb{N}_{a} =\lbrace a, a+1,...\rbrace$ and $ _{b}\mathbb{N} =\lbrace b,b-1,...\rbrace.$ Throughout this paper,  we apply
   discrete nabla and delta  distributions   having supports $\mathbb{N}_{1} $ and $\mathbb{N}_{\alpha-1},\, \alpha>0, $ respectively. 
    For
further details about these distributions, we refer the reader to Ganji and Gharari (2018a). The  dalta and nabla discrete distributions include an  extended class of discrete distributions. The geometric distribution (the number of independent trials required for first success), nabla discrete gamma distribution (Ganji and Gharari (2018a)), nabla  discrete Weibull distribution (Ganji and Gharari (2016)), and nabla  Mittag-Leffler distribution (Ganji and Gharari (2018b)) are examples of nabla discrete distributions. Also, examples for delta discrete distribution are the geometric distribution (the number of failures  for first success), delta discrete gamma distribution, delta discrete Weibull distribution, and delta Mittag-Leffler distribution.
\\
Suppose that $ X $ is a delta discrete random variable  with values $ x= \mathbb{N}_{\alpha-1},\,\,  \alpha>0. $ The delta discrete Laplace transform of probability mass function (pmf) of $ X $ is defined as 
\begin{equation*}
\mathcal{L}_{X}(t)=\sum_{x=\alpha-1} ^{\infty} (\frac{1}{1+t})^{\sigma (x)} f(x)=E[(1+t)^{-\sigma(X)}]=E[e_{\ominus t}(\sigma(X),0)]:=M_{\sigma(X)}(-t),
\end{equation*}
where  $ t $ is the set of all regressive complex constants for which the series converges and $\sigma(x)=x+1.$ By appling the series expansion for  $ (1+t)^{-x}, $  it can be easily 
proved that $ M_{\sigma(X)}(-t)=\sum_{k=0}^{\infty}E[(-1)^{k}(\sigma(X))^{^{\frac{}{k}}}]\frac{t^{k}}{k!},$ 
 where 
\begin{equation*}
(\sigma(X))^{^{\frac{}{k}}}=\dfrac{\Gamma(\sigma(X)+k)}{\Gamma(\sigma(X))}.
\end{equation*}
This function generates the nabla moments of integer order of $ X $  as
\begin{equation*}
 E[(\sigma(X))^{^{\frac{•}{K}}}]=(-1)^{k} 
\frac{d^{k}M_{\sigma(X)}(-t)}{dt^{k}}\mid_{t=0}.
 \end{equation*}
\begin{definition}
Let $ X $ be a delta discrete random variable with pmf $  f.$ \\
$ (a) $ Its $ k$th nabla moment  is denoted by $ \mu_{k}^{\nabla} $ and is defined by $\mu_{k}^{\nabla}:=\sum_{x} (\sigma(X))^{^{\frac{•}{k}}} f(x).$\\
$ (b) $ The nabla moment generating function (mgf)  of $ X $ is given by
 $ M_{\sigma(X)}(t) =E[e_{ t}^{\ast}(\sigma(X),0)],$
 where $e_{ t}^{\ast}(\sigma(X),0)$ denotes the  nabla  exponential function  defined on a time scale (see Bohner and Peterson (2003)). 
\end{definition}
\begin{definition}\label{its}
Let $X$ and $Y$ be delta discrete random variables with pmfs  $f $ and $g, $  distribution functions (dfs) $F$ and $G$ and survival functions (sfs) $\bar{F}$ and $\bar{G},$ respectively. Then $ X$ is said to be smaller than $Y$ in:\\
$(a) $ delta Laplace transform order (written as $X \leq_{Lt}Y $) if $\mathcal{L}_{X}(s)\geq \mathcal{L}_{Y}(s),\,\,\,\,\forall s>0;$\\
$(b) $ delta Laplace transform ratio order (written as $X\leq_{Lt-r}Y$) if $ \mathcal{L}_{X}(s) / \mathcal{L}_{Y}(s) $ is increasing  in $ s>0;$\\
$(c) $ reverse delta Laplace transform ratio order (written as $X \leq_{r-Lt-r}Y$) if $ (1-\mathcal{L}_{X}(s)) / (1-\mathcal{L}_{Y}(s)) $ is increasing  in $ s>0;$\\
$(d) $ \,$i-$differentiated delta Laplace transform ratio order (I) (written as $X \leq_{d^{^{(i)}}-Lt-r}Y$), $i>0,$ if $ \mathcal{ L }^{^{(i)}}_{Y}(s) / \mathcal{ L }^{^{(i)}}_{X}(s) $ is decreasing  in $ s>0.$
\end{definition}
Now, suppose that $ X $ is a nabla discrete random variable with values $ x= \mathbb{N}_{1}.$  The nabla discrete Laplace transform of pmf of $ X $ is defined as 
\begin{equation*}
\mathcal{L}_{X}(t)=\sum_{x=1}^{\infty}(1-t)^{\rho(x)}f(x)=E[(1-t)^{\rho(x)}]= E[e_{\ominus t}^{\ast}(\rho(X),0)]:=M_{\rho(X)}(-t),
\end{equation*}
where $\rho(x)=x-1$ and $ t $ is the set of all regressive complex constants, for which  the series converges. By using the series expansion for  $ (1-t)^{x} ,$  it can be easily 
proved that $ M_{\rho(X)}(-t)=\sum_{k=0}^{\infty} E[(-1)^{k}(\rho(X))^{^{\frac{k}{}}}]\frac{t^{k}}{k!}, $ where
\begin{equation*}
(\rho(X))^{^{\frac{k}{•}}}=\dfrac{\Gamma(\rho(X)+1)}{\Gamma(\rho(X)+1-k)}.
\end{equation*}
This function generates the delta moments of integer order of $ X $  as 
\begin{equation*}
E[(\rho(X))^{^{\frac{K}{}}}]=(-1)^{k} 
\frac{d^{k}M_{\rho(X)}(-t)}{dt^{k}}\mid_{t=0}. 
\end{equation*}
\begin{definition}
Let $ X $ be a nabla discrete random variable with pmf
 $  f.$ \\
$ (a) $ Its $ k $th delta moment  is denoted by $ \mu_{k}^{\Delta} $ and is defined by $\mu_{k}^{\Delta}:=\sum_{x} (\rho(X))^{^{\frac{k}{}}} f(x).$\\
$ (b) $ The delta  mgf  of $ X $ is given by
 $ M_{\rho(X)}(t) =E[e_{ t}(\rho(X),0)],$
  where $e_{ t}(\rho(X),0)$  denotes the  delta  exponential function  defined on a time scale (see Bohner and Peterson (2001, 2003)).
\end{definition}
Also, $X\leq_{Lt}Y,\,X\leq_{Lt-r}Y,\,X\leq_{r-Lt-r}Y$ and $X \leq_{d^{^{(i)}}-Lt-r}Y$  are defined for nabla  random variables $ X $ and $ Y $ like  delta  random variables, but in this case there are two differences;  these results are valid for the restriction $0 <s<1, $  and  Definition \ref{its} (d)  is valid for $i \in\mathbb{N}.$\\
It would be interesting to see a   unification of the  mgfs and  moments.
For a given  time scale $ \mathbb{T}, $ we present the construction  of moments and the mgfs on time scales as 
\begin{equation}\label{cco}
 \hat{\mu}_{k}=E[\Gamma(k+1) \hat{h}_{k}(\eta(X))] \,\,\,\,and\,\,\,\,M_{\eta(X)}(-t)=E[\hat{e}_{\ominus t}(\eta(X),0)],\,\,\,\,\,\,x\in \mathbb{T},
\end{equation}
 respectively. In order that  the reader sees how ordinary moments, delta and nabla moments and also types of the mgfs  follow from  (\ref{cco}),  it is necessary only  in this point to know that 
\begin{equation*}
\hat{h}_{k}(x)=h_{k}(x)=\frac{x^{k}}{\Gamma(k+1)},\,\,\,\eta(x)=\sigma(x)=\rho(x)=x\,\,\,and \,\,\,\hat{e}_{\ominus t}(\eta(X),0)=e^{-tx},
\end{equation*}
if $ \,\,\mathbb{T}=\mathbb{R^{+}}, $ 
\begin{equation*}
\hat{h}_{k}(x)=h_{\frac{k}{•}}(x)=\frac{x^{\frac{k}{•}}}{\Gamma(k+1)},\,\,\,\,\,\eta(x)=\sigma(x)\,\,\,\,\,and\,\,\,\,\,\hat{e}_{\ominus t}(\eta(X),0)=(1+t)^{-\sigma(x)},
\end{equation*}
if $ \,\,\mathbb{T}=\mathbb{N}_{\alpha-1},\,\,  \alpha>0 $  and
\begin{equation*}
\hat{h}_{k}(x)=h_{\frac{}{k}}(x)=\frac{x^{\frac{}{k}}}{\Gamma(k+1)},\,\,\,\,\,\eta(x)=\rho(x)\,\,\,\,\,\,\,and\,\,\,\,\,\,\,\hat{e}_{\ominus t}(\eta(X),0)=(1-t)^{\rho(x)},
\end{equation*}
if $ \,\,\mathbb{T}=\mathbb{N}_{1}. $ 
\begin{definition}\label{kim}
 It is said that the random variable $ X $ has a delta discrete gamma distribution with $ (\alpha,\beta) $ parameters if its pmf is given by
\begin{align}\label{ee}
Pr[X=x]=\frac{h_{\frac{\alpha-1}{•}}(x)\beta^{\alpha}}{e_{\beta}(\sigma(x),0)}=\frac{x^{\frac{\alpha-1}{•}}\beta^{\alpha}}{\Gamma(\alpha)(1+\beta)^{\sigma(x)}},\,\,\,\,\,\,\,\,\,x=\mathbb{N}_{\alpha-1},
\end{align}
where $ \alpha>0,\,\, \beta >0 $ and  it is denoted by $ \Gamma^{\Delta} (\alpha,\beta).$
\end{definition}
\begin{definition}\label{bibi}
It is said that the random variable $ X $ has a nabla discrete gamma distribution with $ (\alpha,\beta) $ parameters if its pmf is given by
\begin{align}\label{hh}
Pr[X=x]=\frac{h_{\frac{}{\alpha-1}}(x)\beta^{\alpha}}{e_{\beta}^{\ast}(\rho(x),0)}=\frac{x^{\frac{}{\alpha-1}}\beta^{\alpha}(1-\beta)^{\rho(x)}}{\Gamma(\alpha)},\,\,\,\,\,\,\,\,\,x=\mathbb{N}_{1},
\end{align}
where $ \alpha>0,\,\,0 <\beta <1 $ and  it is denoted by $ \Gamma^{\nabla} (\alpha,\beta).$
\end{definition}
\section{A new definition of fos by spacings}
Suppose that $n-1$ points are dropped at random on the unit  interval $ (0,1).$ 	We denote the ordered distances of these points from the origin by $u_{i:n}\,(i=1,2,...,n-1).$ We let $V_{i}:=\nabla u_{i:n}=u_{i:n}-u_{i-1:n},\,(u_{0:n}=0) $ and $V^{^{'}}_{i}:=\triangle u_{i:n}=u_{i+1:n}-u_{i:n} $ and call  \textit{backward spacing}  and \textit{forward spacing}, respectively. 
Similarly,
we  may have the same results for $V_{i:n}^{^{,}}\,(i=1,2,...,n-1)$   as
\begin{equation*}
 W_{i}:=\nabla V_{i:n}=V_{i:n}-V_{i-1:n}=u_{i:n}-2u_{i-1:n}+u_{i-2:n}=\nabla^{2}u_{i:n}
\end{equation*}
and
\begin{equation*}
 W^{^{'}}_{i}:=\triangle V_{i:n}=V_{i+1:n}-V_{i:n}=u_{i+1:n}-2u_{i:n}+u_{i-1:n}=\triangle^{2}u_{i:n}.
\end{equation*}
By repeating $m$ times, we obtain
\begin{equation*}
\nabla^{m}u_{i:n}=\sum _{j=0}^{m}(-1)^{j}\binom{m}{j}u_{i-m-j:n}:=Z_{i}
\end{equation*}
and
\begin{equation*}
\triangle ^{m}u_{i:n}=\sum _{j=0}^{m}(-1)^{j}\binom{m}{j}u_{i+m-j:n}:=Z^{^{'}}_{i},
\end{equation*}
where $\nabla^{m}$ and $\triangle^{m}$ are $m$th order nabla and delta  difference operators, respectively.
\\
Now let $\alpha>0$ with $n^{'}-1<\alpha<n^{'}. $ using the definition of nabla and delta differences of order $\alpha $ from discrete fractional calculus (see Goodrich and Peterson (2015)), we have
\begin{equation*}
\nabla^{\alpha}u_{i:n}=\sum _{j=0}^{i-2}(-1)^{j}\binom{\alpha}{j}u_{i-j:n}:=S_{i},\,\,\,\,\,\,i\in\mathbb{N}_{1+n^{'}}
\end{equation*}
and
\begin{equation*}
\triangle ^{\alpha}u_{i:n}=\sum _{j=0}^{\alpha+i-1}(-1)^{j}\binom{\alpha}{j}u_{i+\alpha-j:n}:=S^{^{'}}_{i},\,\,\,\,\,\,i\in\mathbb{N}_{1+n^{'}-\alpha}.
\end{equation*}
Considering the indices of $S^{^{'}}_{i} $ and $S_{i}, $
   we proved   new definitions of fos and os, respectively.
Therefore, we could build  a new definition of fos  by using discrete fractional calculus (or nabla and delta fractional calculus). On the other hand, the existence of the delta and nabla discrete distributions were proved by using discrete fractional calculus (see Ganji and Gharari (2018)). It seems that there is a relationship between os and fos with these types of probability distributions. In the following, we are going to discover it with an interesting way.
\section{Distribution of a fos by fractional operators }
 Stigler (1977) considered fos with both integer and non-integer sample sizes. Rohatgi and Saleh (1988) extended definition of df of the $ i$th os, $ i \in \mathbb{N}, $ with noninteger sample size. In this section, first we give two general expressions for the pdfs of the $ i$th os and fos
when the sample size $ N $ is a  nabla  or delta  random variable.  For doing this, we use both  ordinary  and fractional types of $i$th derivative of  Laplace transform of $ N.$ 
   Then,
  we introduce a new stochastic
ordering based on the $i$th derivative 
 of the Laplace transform $ N.$
\subsection{Distribution of an os when the  sample size is random }
 Let $X_{1},X_{2},...,X_{N}$ be a random sample  size $ N $ from a df  $ F. $  Let  $X_{1:N},X_{2:N},...,X_{N:N}$  be the corresponding os when the sample size $N$ is a random variable itself. 
 \\
Conditionally, given sample size $N =k,$ let the df
of the $i$th os $X_{i:k}$ be noted by $ F_{i:k}(x \vert k)$ and the corresponding pdf by $f_{i:k}(x \vert k).$
We suppose that $ N $ is a random variable with a nabla discrete distribution. In this case, 
\begin{equation*}
P_{N}(k)=P(N=k),\,\,\,\,\,\,k\in\mathbb{N}_{1}.
\end{equation*}
The conditional df of $X_{i: \rho(N) },$  conditioned
on the event $ \lbrace N \geq i \rbrace, $  is 
\begin{equation}\label{b}
F_{i:\rho(N)}(x )=\frac{1}{\pi_{N}^{i}}\sum_{k=i}^{\infty}F_{i:\rho(k)}(x\vert k)P_{N}(k),
\end{equation}
where $\pi_{N}^{i}=P(N\geqslant i).$
 It is well known that  the df of the  $i$th os, $i\in\mathbb{N}_{1}, $ in a sample  size  $ \rho(k) $ from a  df  $F $ is as
\begin{equation}\label{ufo1}
F_{i:\rho(k)}(x\vert k)=\sum_{j=i}^{\rho(k)}\binom{\rho(k)}{j}F^{j}(x)(1-F(x))^{\rho(k)-j}.
\end{equation}
Obviously,
\begin{equation*}
F_{i:\rho(k)}(x\vert k)=\int_{0}^{F(x)}\dfrac{\rho(k)!}{(i-1)!(\rho(k)-i)!}t^{i-1}(1-t)^{\rho(k)-i}dt.
\end{equation*}
 The $i$th derivative of $M_{\rho(N)}(t)$ exists and is obtained by formal differentiation as
\begin{equation}\label{a}
(-1)^{i}M_{\rho(N)}^{(i)}(t)=\sum_{k=i}^{\infty}(\rho(k))^{^{\frac{i}{•}}}(1-t)^{\rho(k)-i}P_{N}(k),\,\,\,\,\,i\in\mathbb{N}_{1}.
\end{equation}
It follows from (\ref{a}) and (\ref{b}) that for $x\in\mathbb{R},$
\begin{align*}
F_{i:\rho(N)}(x)
&=\frac{1}{\pi_{N}^{i}}\sum_{k=i}^{\infty}\int_{0}^{F(x)}\frac{(\rho(k))^{^{\frac{i}{•}}}}{\Gamma(i)}t^{i-1}(1-t)^{\rho(k)-i}P_{N}(k)dt\\
&=\frac{1}{\pi_{N}^{i}(i-1)!}\int_{0}^{F(x)}(-1)^{i}t^{i-1}M_{\rho(k)}^{(i)}(t)dt,\,\,\,\,\,i\in\mathbb{N}_{1}.
\end{align*}
In the special case, when $ F$ has pdf $ f, $ we can immediately  get the pdf of the $ i$th os when the  sample size  is random as
\begin{equation}\label{c}
f_{i:\rho(N)}(x)=\frac{(-1)^{i}F^{i-1}(x)f(x)}{\pi_{N}^{i}(i-1)!}M_{\rho(N)}^{(i)}(F(x)),\,\,\,\,\,x\in\mathbb{R}.
\end{equation}
Note that this formula is valid for $ i=1,2,..., \rho(N).$ Now, let us introduce a new order  which, as will be seen in the next section, has some potential applications.
\begin{definition}
Let $ X $ and $Y$ be nabla  random variables with pdfs  $f $ and $g, $ dfs $F$ and $G,$ and sfs $\bar{F}$ and $\bar{G},$ respectively. Then, $ X$ is said to be smaller than $Y$ in the
 $ i-$differentiated nabla Laplace transform ratio order (II)  (written as $ X\leq_{ D ^{^{i}}-Lt-r}Y $) if\\

$\,\,\,\,\,\,\,\,\,\,\,\,\,\,\,\,\,\,\,\,\,\,\,\,\,\,\,\,\,\,\,\,\,\,\,\,\,\,\,\,\,\,\,\,\,\,\,\,\,\,\,\,\,
\dfrac{\pi _{X}^{i}\mathcal{L}^{(i)}_{Y}(s)}{\pi _{Y}^{i}\mathcal{L}^{(i)}_{X}(s)} \,\,$ is decreasing  in $ 0<s<1,\,\,i\in\mathbb{N}_{1}. $ \\
In the special case $i=1, $ we have $P(X\geqslant 1)=P(Y\geqslant 1)=1$ and  the order reduces to  $\,\,\dfrac{\mathcal{\acute{L}}_{Y}(s)}{\mathcal{\acute{L}}_{X}(s)} $ (written as $X \leq_{d-Lt-r}Y$). 
\end{definition}

\subsection{Distribution of a fos when the  sample size is fractional  random}
Let $X_{1},X_{2},...,X_{N}$ be a random sample  size $ N $ from a df  $ F $ with  corresponding sf $\bar{F}. $
Suppose that $ N $ is a random variable with a delta discrete distribution. In this case,
\begin{equation*}
P_{N}(k)=P(N=k),\,\,\,\,\,\,k\in\mathbb{N}_{\alpha-1},\,\,\,\,\alpha>1,
\end{equation*}
is the distribution of $N.$
\\
We let  $X_{\gamma:\sigma(N)+\gamma-n},\,\, n-1<\gamma\leq n,\,n\in\mathbb{N}$ be the corresponding fos when     $ N$   is  a delta random variable itself. Note that, in the special cases $\gamma=1, 2,...,n\,, $ we get $X_{1:\sigma(N)},\,X_{2:\sigma(N)},\,..., \,X_{n:\sigma(N)}$ os,  respectively. Also, in the more special case, if  $\alpha=2$ and the sample size is constant $n,$   for   $\gamma=1,\,2,\,...\,n\,, $ we get $ X_{1:n},\,X_{2:n},\,...,\,X_{n:n}$ os,  respectively.
\\We are interested in obtaining the df of the $\gamma$th fos  in the fractional sample  size  $ \sigma(k)+\gamma-n,\,\,n-1<\gamma\leq n, $ from a  df  $F $  when   $N$ has a  delta discrete  distribution. Also, when $F$ has pdf $f, $ an expression  for the pdf $f_{\gamma:\sigma(N)+\gamma-n}$ of the $\gamma$th fos is derived. 
Conditionally, given sample size with $N =k,$ the df
of the $\gamma$th fos $X_{\gamma:\sigma(k)+\gamma-n}$ is noted by $ F_{\gamma:\sigma(k)+\gamma-n}(x \vert k )$ and the corresponding pdf by $f_{\gamma:\sigma(k)+\gamma-n}(x \vert k ),$ where $ n-1<\gamma\leq n.$
We write  the df of the $\gamma$th fos  in a fractional sample  size  $ \sigma(k)+\gamma-n $ from a  df  $F $ as
\begin{equation}
F_{\gamma:\sigma(k)+\gamma-n}(x\vert k)=\sum_{j=\gamma}^{\infty}\binom{\sigma(k)+\gamma-n}{j}F^{j}(x)(1-F(x))^{ (\sigma(k)+\gamma-n)-j}.
\end{equation} 
This is actually the same formula as before (\ref{ufo1}), since $\binom{\sigma(k)+\gamma-n}{j}=0 $  for $\sigma(k)+\gamma-n<j$ when $j$ is a integer or $\gamma=n.$ As an example, for $\gamma=1,$ we have the df of minimum os as
\begin{equation*}
F_{1:\sigma(k) }(x\vert k)=\sum_{j=1}^{\sigma(k) }\binom{\sigma(k) }{j}F^{j}(x)(1-F(x))^{  \sigma(k) -j}.
\end{equation*}
 Obviously,   
\begin{equation*}
F_{\gamma:\sigma(k)+\gamma-n}(x\vert k)=\int_{0}^{F(x)}\dfrac{\Gamma(\sigma(k)+\gamma-n+1)}{\Gamma(\gamma )\Gamma(\sigma(k)-n+1)}t^{\gamma-1}(1-t)^{\sigma(k)-n}dt.
\end{equation*}
 Now,
 we obtain
 fractional derivative of order $\gamma$ of $M_{\sigma(N)}(t).$  For further details on fractional operators, we refer the reader to Miller and Ross (1993). Because $f(t)=(1+t)^{-\sigma(k)}$ is a function having a power series representation
 $\sum_{i=0}^{\infty}\dfrac{(\sigma(k))^{^{\frac{•}{i}}}(-1)^{i}t^{i}}{i!}$
and also
$\mathcal{D}^{\gamma} _{t} t^{i}=\dfrac{\Gamma(i+1)}{\Gamma(i+1-\gamma)}t^{i-\gamma}, $ where $\mathcal{D}^{\gamma} _{t}$ denotes the $\gamma$th fractional derivative, we can write the $\gamma$th fractional derivative of this function as
  \begin{equation*}  
 \mathcal{D}^{\gamma} _{t} (1+t)^{-\sigma(k)}=\sum_{i=\gamma}^{\infty}\dfrac{(\sigma(k))^{^{\frac{•}{i}}}(-1)^{i}t^{i-\gamma}}{\Gamma(i+1-\gamma)}.
\end{equation*} 
On the other hand, we have
 \begin{equation*}  
  \sum_{i=\gamma}^{\infty}\dfrac{\Gamma(k+1+i)(-1)^{i}t^{i-\gamma}}{\Gamma(k+1)\Gamma(i+1-\gamma)}=(-1)^{\gamma}(\sigma(k))^{^{\frac{•}{\gamma}}}(1+t)^{-\sigma(k)-\gamma}.
\end{equation*} 
 So, we obtain
\begin{equation}\label{ufo2}
(-1)^{\gamma}\mathcal{D}^{\gamma} _{t} M_{\sigma(N)}(t)=\sum_{k=\gamma}^{\infty}(\sigma(k))^{^{\frac{\gamma}{•}}}(1+t)^{-\sigma(k)-\gamma }P_{N}(k),\,\,\,\,\,k\in\mathbb{N}_{\alpha-1}.
\end{equation}
The conditional df of $X_{\gamma: \sigma(N)+\gamma-n  }$  
on the event $ \lbrace N \geq \gamma \rbrace $  is 
\begin{equation}\label{ufo3}
F_{\gamma:\sigma(N)+\gamma-n }(x )=\frac{1}{\pi_{N}^{\gamma}}\sum_{k=\gamma}^{\infty}F_{\gamma:\sigma(k)+\gamma-n}(x\vert k)P_{N}(k),
\end{equation}
where $\pi_{N}^{\gamma}=P(N\geqslant \gamma).$
It follows from (\ref{ufo2}) and (\ref{ufo3}) that for $x\in\mathbb{R},$
\begin{align*}
F_{\gamma:\sigma(N)+\gamma- n}(x)
&=\frac{1}{\pi_{N}^{\gamma}}\sum_{k=\gamma}^{\infty}\int_{0}^{F(x)}\frac{\Gamma(\sigma(k)+\gamma-n+1)}{\Gamma(\gamma)\Gamma(\sigma(k)-n+1)}t^{\gamma-1}(1-t)^{\sigma(k)-n}P_{N}(k)dt\\
&=\frac{1}{\pi_{N}^{\gamma}\Gamma(\gamma) }\int_{0}^{\frac{F }{\bar{F}}}\frac{\Gamma(\sigma(k)+\gamma-n+1)}{\Gamma(\gamma)\Gamma(\sigma(k)-n+1)}t^{\gamma-1}(1+t)^{-(\sigma(k)-n+1)-\gamma}P_{N}(k)dt\\
&=\frac{1}{\Gamma(\gamma)\pi_{N}^{\gamma}}\int_{0}^{\frac{F }{\bar{F}}}
(-1)^{\gamma}t^{\gamma-1} \mathcal{D}^{\gamma} _{t} M_{\sigma(\sigma(k)-n)}(t)dt,\,\,\,\,\,\gamma\in\mathbb{N}_{\alpha-1}.
\end{align*}
In the special case, when $ F$ has pdf $ f, $ we can immediately  get the pdf of the $ \gamma$th fos as
\begin{equation}\label{c2}
f_{\gamma:\sigma(k)+\gamma-n}(x)=\frac{(-1)^{\gamma}f(x)F^{\gamma-1}(x)}{\pi_{N}^{\gamma}\Gamma(\gamma)\bar{F}^{\gamma+1}(x)} \mathcal{D}^{\gamma} _{t} M_{\sigma(\sigma(N)-n)}(\frac{F(x)}{\bar{F}(x)}),\,\,\,\,\,x\in\mathbb{R}.
\end{equation}
Note that in the special case $\gamma=1/2$ we get $(-1)^{\gamma}=i, $ where $i$ is   unit imaginary number.
 In the following, we are interested in  comparing  fos when the sample size is a random variable. Toward this end,
 we define a new stochastic order based on the $\gamma $th fractional derivative  of the Laplace transform of $ N. $ 

\begin{definition}
Let $X$ and $Y$ be delta discrete random variables with pdfs  $f $ and $g, $ dfs $F$ and $G,$ and sfs $\bar{F}$ and $\bar{G},$ respectively. Then, $ X$ is said to be smaller than $Y$ in the
 $\gamma-$differentiated delta Laplace transform ratio order (II) (written as $ X\leq_{D^{^{\gamma}}-Lr-r}Y $) if\\
 
$\,\,\,\,\,\,\,\,\,\,\,\,\,\,\,\,\,\,\,\,\,\,\,\,\,\,\,\,\,\,\,\,\,\,\,\,\,\,\,\,\,\,\,\,\,\,\,\,\,\,\,\,\,
\dfrac{\pi _{X}^{\gamma}\mathcal{L}^{(\gamma)}_{Y}(s)}{\pi _{Y}^{\gamma}\mathcal{L}^{(\gamma)}_{X}(s)}\,\,\,$ is decreasing  in $ s>0,\,\,\gamma\in\mathbb{N}_{\alpha-1},\,\,\alpha>1. $ 
\end{definition}
In the special case $\gamma=\alpha-1, $ we have $P(X\geqslant \alpha-1)=P(Y\geqslant \alpha-1)=1$ and   the recent definition  reduces to $ \mathcal{L}^{(\alpha-1)}_{Y}(s)/\mathcal{L}^{(\alpha-1)}_{X}(s). $
\\
For the rest of this paper, we  denote the $i$th os with $X_{i:N}$ when the sample size $N$ is a nabla random variable,  and the  $\gamma$th fos  with $X_{\gamma:N}$  when $N$ is a delta random variable  in the fractional  sample size $\sigma(N)+\gamma-n.$  

\section{A new stochastic order based on discrete Laplace transform }

If  $ X $ is a delta random variable with values $\mathbb{N}_{\alpha-1},$ we denote its delta Laplace transform with $\mathcal{L}_{X}(s) $ and we have
 \begin{equation*} 
 \Psi_{X}(s) = 1-\mathcal{L}_{X}(s) =\int^{\infty}_{\alpha-1}(1-e_{\ominus s}(\sigma(t),0))f(t)\Delta t
\end{equation*} 
 or
  \begin{equation*} 
 \Psi_{X}(s) = \sum^{\infty}_{t=\alpha-1}(1-(\frac{1}{1+s})^{\sigma(t)})f(t),
\end{equation*} 
which can be considered as a mixture distribution of the Burr  distribution
with mean $ 1/t,\,t>-1 $ 
 and mixing distribution $ F. $ Thus, $\Psi$ is a distribution with a
density
\begin{equation*} 
 \psi_{X}(s) =\int _{\alpha-1}^{\infty}\sigma(t)(1+s)^{\sigma(t)-1}f(t)\Delta t,\,\,\,\, s>0,\,t>-1.
 \end{equation*} 
  The random variable with this distribution  is denoted by $\xi(X). $ Assume that $X $ and $Y $ are delta discrete random variables with 
 density functions $f $ and $g, $ respectively. Recall that $X $ is smaller than $Y $ in the likelihood ratio order (denoted by $ X \leq _{lr} Y$)
if $g(x)/f (x) $ is increasing in $ x. $ Note that $-\mathcal{\acute{L}}_{X}(s)$ and $-\mathcal{\acute{L}}_{Y}(s)$
are density functions  of $\xi(X) $ and $\xi(Y), $ respectively. It is easy to see  that
\begin{align}\label{won2}
X\leq _{d-Lt-r} Y \Leftrightarrow \xi (Y) \leq _{lr} \xi (X)
\end{align}
and
\begin{align}\label{won1}
X \leq _{Lt (Lt-r, r-Lt-r)} Y\Leftrightarrow \xi (Y)\leq _{st(hr,rh)}\xi (X).
\end{align}
The definitions of the stochastic orders that are mentioned above can be found, for example, in Shaked and Shanthikumar (1994).
 \\
Consider the Laplace transform of $F$ as
\begin{equation*} 
\mathcal{L}^{*}_{X}(s)=\int^{\infty}_{\alpha-1}e_{\ominus s}(\sigma(t),0)F(t)\Delta t
\end{equation*}
and define the Laplace transform  of $\bar{F}$ as
\begin{equation*} 
\mathcal{L}^{**}_{X}(s)=\int^{\infty}_{\alpha-1}e_{\ominus s}(\sigma(t),0)\bar{F}(t)\Delta t,
\end{equation*}
for all $ s>0. $ By applying Corollary 2.14 and Example 2.5 in Goodrich and Peterson (2015), it is easy to verify 
\begin{align} \label{hi}
\mathcal{L}_{X}^{*}=\frac{\mathcal{L}_{X}}{s(1+s)^{\alpha-1}},\,\,\, \mathcal{L}_{X}^{**}=\frac{1-\mathcal{L}_{X}}{s(1+s)^{\alpha-1}}.
\end{align}
 Also, $\mathcal{L}_{Y},\,\mathcal{L}_{Y}^{*} $ and $\,\mathcal{L}_{Y}^{**}$ are defined for $Y$ like $X.$\\
Let $X$ be a delta random variable with df $F$
and sf $\bar{F} = 1 - F.$
Denote $R_{s}^{X}(n)=(-1)^{n-1}s^{n}/(n-1)!(d^{n-1}/ds^{n-1})\mathcal{L}^{**}_{X}(s),\, n\geq 1,\, 0<s,$ and let $R_{s}^{X}(0)=1$ for all $0<s.\,\,R_{s}^{X}(n) $ is the reliability function corresponding to a discrete random variable, say, $ N_{s}^{X}.$ Note that, by Theorem 1.58 (Integration by part) in Goodrich and Peterson (2015),  $R_{s}^{X}(n)$
 can be written as
\begin{align*}
R_{s}^{X}(n)=&\int_{\alpha-1}^{\infty}\sum_{i=n}^{\infty}(1-\frac{s}{1+s})^{\sigma(x)}\frac{(\sigma(x))^{\frac{}{i}}}{i!}(s/1+s)^{i}\Delta F(x)\\
&=\left\{
{\begin{array}{*{20}c}
\int_{\alpha-1}^{\infty}\Gamma^{\nabla}_{\eta}(n,\sigma(x) )\bar{F}( x)\Delta x,\,\,\,\,\,\, \,\,\,\,\,\,\,\,\,\,\,\,\,&n=1,2,...\\
\,\,\,\\
1,\,\,\,\,\,\,\,\,\,\,\, \,\,\,\,\,\,\,\,\,& n=0\\
\end{array}}\right.
\end{align*}
where $\eta=s/1+s.$\\
The following result gives a delta Laplace transform characterization of the order $\leq_{st}.$
\begin{theorem}
Let $X $ and $Y$ be two delta random variables, and
let $N_{s}(X )$ and $N_{s}(Y)$ be as described above. Then,
\begin{equation*}
X \leq_{st} Y\Rightarrow N_{s}(X )\leq_{st}N_{s}(Y).
\end{equation*}
\end{theorem}
If  $ X $ is a nabla random variable with values $\mathbb{N}_{1},$ we denote its nabla Laplace transform with $\mathcal{L}_{X}(s) $ and we have
 \begin{equation*} 
 \Psi_{X}(s) = 1-\mathcal{L}_{X}(s) =\int^{\infty}_{0}(1-e_{\ominus s}^{*}(\rho(t),0))f(t)\nabla t
\end{equation*} 
 or
  \begin{equation*} 
 \Psi_{X}(s) = \sum^{\infty}_{t=1}(1-(1-s)^{t-1})f(t),
\end{equation*} 
which can be considered as a mixture distribution of the beta distribution with mean $ 1/t, (t>1) $ and a mixing distribution $ F. $ Thus, $\Psi$ is a distribution with a
density  funtion   
\begin{equation*} 
 \psi_{X}(s) =\int _{0}^{\infty}\rho(t)(1-s)^{\rho(t)-1}f(t)\nabla t,\,\,\,\,0<s<1, \,t>1.
 \end{equation*} 
 The random variable with this distribution is denoted by $\xi(X).$  Obviously, $-\mathcal{\acute{L}}_{X}(s)$ and $-\mathcal{\acute{L}}_{Y}(s)$
are density functions  of $\xi(X) $ and $\xi(Y), $ respectively. Also, the relations (\ref{won1}) and (\ref{won2}) are valid for  random variables  $\xi(X) $ and $\xi(Y). $
Consider the Laplace transform of $F$ as
\begin{equation*} 
\mathcal{L}_{X}^{*}(s) =\int^{\infty}_{0}e_{\ominus s}(\rho(t),0)F(t)\nabla t
\end{equation*}
and define the Laplace transform  of $\bar{F}$ as
\begin{equation*} 
\mathcal{L}_{X}^{**}(s) =\int^{\infty}_{0}e_{\ominus s}(\rho(t),0)\bar{F}(t)\nabla t,
\end{equation*}
for all $ 0<s<1. $ By using Theorem 3.82 and Example 3.36 in Goodrich and Peterson (2015), it is easy to verify 
\begin{equation}\label{labb}
\mathcal{L}_{X}^{*} =\frac{1}{s}\mathcal{L}_{X} ,\,\,\, \mathcal{L}_{X} ^{**}=\frac{1-\mathcal{L}_{X} }{s}. 
\end{equation}
Also, $\mathcal{L}_{Y},\,\mathcal{L}_{Y}^{*} $  and $\,\mathcal{L}_{Y}^{**}$ are defined for $Y$ like $X.$\\
Let $X$ be a nabla random variable with distribution function $F$
and sf $\bar{F}= 1 - F.$
Denote $R_{s}^{X}(n)=(-1)^{n-1}s^{n}/(n-1)!(d^{n-1}/ds^{n-1})\mathcal{L}^{**}_{X}(s),\,\,n\geq 1,\,\,0<s<1,$ and let $R_{s}^{X}(0)=1$ for all $0<s<1.\,\,R_{s}^{X}(n) $ is the reliability function corresponding to a discrete random variable, say, $ N_{s}^{X}.$ Note that, by Theorem 3.36 (Integration by part) in \cite{15}, $R_{s}^{X}(n)$
 can be written as
\begin{align*}
R_{s}^{X}(n)=&\int_{0}^{\infty}\sum_{i=n}^{\infty}(1+\frac{s}{1-s})^{-\rho(x)}\frac{(\rho(x))^{\frac{i}{•}}}{i!}(s/1-s)^{i}\nabla F(x)\\
&=\left\{
{\begin{array}{*{20}c}
  \int_{0}^{\infty}\Gamma^{\Delta}_{\eta}(n,\rho(x) )\bar{F}( x)\nabla x,\,\,\,\,\,\,\,\,\,\,\,\,\,\,\,\, &n=1,2,...\\
\,\\
 1,\,\,\,\,\,\,\,\,\,\,\,\,\,\,\,\,\,\,\,\,\,  & n=0\\
\end{array}}\right.
\end{align*}
where $\eta=s/1-s.$
\\
The following result gives a nabla Laplace transform characterization of the order $\leq_{st}.$
\begin{theorem}
Let $X $ and $Y$ be two nabla random variables, and
let $N_{s}(X )$ and $N_{s}(Y)$ be as described above. Then,
\begin{equation*}
X \leq_{st}Y\Rightarrow N_{s}(X )\leq_{st}N_{s}(Y).
\end{equation*}
\end{theorem}
Using (\ref{hi}) and (\ref{labb}), it is easy to verify the following result.
\begin{theorem}
Let $ X $ and $ Y$ be two  delta (nabla) random variables with
sfs $\bar{F} $and $\bar{G},$ respectively. Then, $X \leq_{Lt} Y$ if  and only if 
$\mathcal{L}_{X}^{**}\leq \mathcal{L}_{Y}^{**}. $
\end{theorem}
If $ X \leq_{Lt} Y, $ then 
\begin{align*}
(1-\mathtt{E} [(1+s)^{\sigma(X)}])/s(1+s)^{\alpha-1}
 &\leq (1 -\mathtt{E}[(1+s)^{\sigma(Y)}])/s(1+s)^{\alpha-1} \\
    (1 - E [ (1 - s)^{ \rho( X )} ] ) / s  &\leq  ( 1 - E [ ( 1 - s) ^ { \rho( Y) } ] ) / s ) 
\end{align*}
 for all
$ s>0\,\, (0<s <1). $ Letting $ s \rightarrow 0, $  it is seen that
\begin{equation}\label{rt}
X \leq_{Lt} Y \Rightarrow \mathtt{E}[\sigma(X)] \leq \mathtt{E}[\sigma(Y)]\,\,\, \left( X \leq_{Lt} Y \Rightarrow \mathtt{E}[\rho(X)] \leq \mathtt{E}[\rho(Y)] \right),  
\end{equation}
provided the expectations exist.
\\
The next proposition characterizes the orders $\leq_{Lt-r}$ and $\leq_{r-Lt-r}$ by functions of the respective moments.
\begin{theorem}
Let $X $ and $Y$ be delta  random variables that possess moments $\mu_{k}^{\nabla}$ and $\nu_{k}^{\nabla},\,(k=1,2,...), $ respectively. Then,\\
 $(a)\,X \leq_{Lt-r}Y$ if and only if

$\,\,\,\,\,\,\,\,\,\,\,\,\,\,\,\,\,\,\,\,\,\,\,\,\,\,\,\,\,\,\,\,\,\,\,\,\,\,\,\,\,\,\,\,\,\,\,\,\,\,\,\,\,\,\,\,\,\,\,\,\,\,\,\,\,\,\,\,\,\,\,\,\,\,\,\,\,\,\,\,\,\,\,\,\,\,\,\,\,\,\,\,\,\,\,\,\dfrac{\sum_{k=0}^{\infty}\frac{(-s)^{k}}{k!}\nu_{k}^{\nabla}}{\sum_{k=0}^{\infty}\frac{(-s)^{k}}{k!}\mu_{k}^{\nabla}}\,\,\,\,$ is decreasing in $ s. $\\
 $(b)\,X \leq_{r-Lt-r}Y$ if and only if

$\,\,\,\,\,\,\,\,\,\,\,\,\,\,\,\,\,\,\,\,\,\,\,\,\,\,\,\,\,\,\,\,\,\,\,\,\,\,\,\,\,\,\,\,\,\,\,\,\,\,\,\,\,\,\,\,\,\,\,\,\,\,\,\,\,\,\,\,\,\,\,\,\,\,\,\,\,\,\,\,\,\,\,\,\,\,\,\,\,\,\,\,\,\,\,\,\,\dfrac{\sum_{k=1}^{\infty}\frac{(-s)^{k}}{k!}\nu_{k}^{\nabla}}{\sum_{k=1}^{\infty}\frac{(-s)^{k}}{k!}\mu_{k}^{\nabla}}\,\,\,\,$ is decreasing in $ s. $
\end{theorem}
\begin{proof}
By using the series expansion $(1-s)^{\rho(t)}=\sum_{k=0}^{\infty}[(-s)^{k}/k!](\rho(t))^{\frac{k}{•}} $ and  $\mu_{0}^{\nabla}=\nu_{0}^{\nabla}=1, $ the result follows easily from Definition \ref{its}.
\end{proof}
\begin{remark}
If  $X $ and $Y$ are nabla  random variables, the same result  holds by substituting $\mu_{0}^{\nabla} $ and $\nu_{0}^{\nabla} $ with
their corresponding nabla types. What we need to  prove is the series expansion $(1+s)^{-\sigma(t)}=\sum_{k=0}^{\infty}[(-s)^{k}/k!](\sigma(t))^{^{\frac{}{k}}} $ and  $\mu_{0}^{\triangle}=\nu_{0}^{\triangle}=1. $
\end{remark}
\begin{theorem}\label{sin}
Let $X$ and  $Y$ be delta (nabla) random variables. If $ X\leq_{Lt-r}Y, $ or if $ X\leq_{r-Lt-r}Y, $ then $\mathcal{L}_{X}(s)\geq\mathcal{L}_{Y}(s)$ for all $ s>0 \, (0<s<1), $  i.e. $ X\leq_{Lt}Y. $
\end{theorem}
\begin{proof}
Denote $\mathcal{L}_{X}(\infty)=\lim_{s \to \infty}\mathcal{L}_{X}(s)\,\,(\mathcal{L}_{X}(1)=\lim_{s \to 1}\mathcal{L}_{X}(s)).$ Since $\mathcal{L}_{X}(0)=1 $ and $\mathcal{L}_{X}(\infty)=0\,\,(\mathcal{L}_{X}(1)=0)$ we see that if $X\leq_{Lt-r}Y, $ then
\begin{equation*}
\frac{\mathcal{L}_{Y}(s)}{\mathcal{L}_{X}(s)}\leq \frac{\mathcal{L}_{Y}(0)}{\mathcal{L}_{X}(0)}=1,
\end{equation*}
and if $X\leq_{r-Lt-r}Y, $ then
\begin{equation*}
\frac{1-\mathcal{L}_{Y}(s)}{1-\mathcal{L}_{X}(s)}\geq \frac{1-\mathcal{L}_{Y}(\infty)}{1-\mathcal{L}_{X}(\infty)}=1\,\,\,\left(\frac{1-\mathcal{L}_{Y}(s)}{1-\mathcal{L}_{X}(s)}\geq \frac{1-\mathcal{L}_{Y}(1)}{1-\mathcal{L}_{X}(1)}=1\right).
\end{equation*}
This proves the stated results.
\end{proof} 
As a corollary of Theorem \ref{sin}, we see that
\begin{align*}
&  X\leq_{Lt-r}Y \rightarrow \mathtt{E}[\sigma(X)] \leq \mathtt{E}[\sigma(Y)],\,\,\,\,  X\leq_{r-Lt-r}Y \rightarrow \mathtt{E}[\sigma(X)] \leq \mathtt{E}[\sigma(Y)]
  \\
 &\left( X\leq_{Lt-r}Y \rightarrow \mathtt{E}[\rho(X)] \leq \mathtt{E}[\rho(Y)],\,\,\,\,   X\leq_{r-Lt-r}Y \rightarrow \mathtt{E}[\rho(X)] \leq \mathtt{E}[\rho(Y)] \right)
\end{align*}
provided the expectations exist; see Equation (\ref{rt}).
It is easily seen that 
\begin{align*}
 -(1-s)\mathcal{\acute{L}}_{X}(s)/\mathtt{E}[\rho(X)],\,\,\,\, \forall 0<s<1\,\,\,\left(-(1+s)\mathcal{\acute{L}}_{X}(s)/\mathtt{E}[\sigma(X)],  \,\,\,\,\forall s>0\right),
 \end{align*} 
 is the nabla (delta) Laplace transform
of the distribution 
\begin{equation*}
\hat{F}(x) =\int _{0}^{x}\rho(u)\nabla F(u)/\mathtt{E}[\rho(X)]\,\,\,\,\left(=\int _{\alpha-1}^{x}\sigma(u)\Delta F(u)/\mathtt{E}[\sigma(X)] \right),
\end{equation*}  
related to $ F, $ provided $ 0<\mathtt{E}[\rho(X)]<\infty\,\,\,\left( 0<\mathtt{E}[\sigma(X)]<\infty\right). $   We denote the random variables corresponding to $\hat{F}$ and $\hat{G}$ by $\hat{X}$ and $\hat{Y},$
 respectively. Then, we get $X\leq_{d-Lt-r}Y$ if and only if $\hat{X}\leq_{Lt-r}\hat{Y}.$ This observation leads to this that $X\leq_{lr}Y$ implies $X\leq_{d-Lt-r}Y. $ Because $X \leq_{lr}Y$ if and only if $ \hat{X}\leq_{lr}\hat{Y}$ and on the other hand $ \hat{X}\leq_{lr}\hat{Y}$ implies
$ \hat{X}\leq_{Lt-r}\hat{Y}. $
\begin{theorem}
Let $n>0$ be a fixed integer. Let $X_{1},\,X_{2},\,...,\,X_{n}$ be a set of independent nabla random variables and let $Y_{1},\,Y_{2},\,...,\,Y_{n}$ be another set of  independent nabla random variables. If $X_{j}\leqslant_{Lt-r}Y_{j},\,j=1,2,...,n, $ then $X_{1}+X_{2}+\,...+\,X_{n}\leqslant_{Lt-r}Y_{1}+Y_{2}+\,...+\,Y_{n}.$
\end{theorem}
\begin{proof}
Since $\mathcal{L}_{X_{1}+X_{2}+\,...+\,X_{n}}(s)=\prod_{i=1}^{n}\mathcal{L}_{X_{i}}(s),$ we see that if $\mathcal{L}_{Y_{j}}(s)/\mathcal{L}_{X_{j}}(s) $ is decreasing in $s,\,j=1,2,...,n, $ then $[\mathcal{L}_{Y_{1}+Y_{2}+\,...+\,Y_{n}}(s)]/[\mathcal{L}_{X_{1}+X_{2}+\,...+\,X_{n}}(s)]$ is also decreasing in $s.$
\end{proof}
 Equations (\ref{c}) and (\ref{c2}) give
the following theorem.  
\begin{theorem}\label{san}
Let $X_{1},\,X_{2},... $ be independent and identically distributed non-negative random variables, and let $N_{1}$ and $N_{2}$ be nabla random variables which are independent of the $X_{i}. $ Then  
\begin{equation*}
N_{1}\leqslant_{Lt-r}\left[\leqslant_{r-Lt-r}\right] N_{2}\Leftrightarrow \sum_{i=1}^{N_{1}}X_{i}\leqslant_{Lt-r} [\leqslant_{r-Lt-r}]\sum_{i=1}^{N_{2}}X_{i}.  
\end{equation*}
\end{theorem}
\begin{proof}
For $j=1,2, $ we have
\begin{align*}
\mathcal{L}_{X_{1}+X_{2}+...+X_{N_{j}}} (s)=&\sum_{k=1}^{\infty}P_{N_{j}}(k)\mathcal{L}_{X_{1}+X_{2}+...+X_{k}} (s)\\
=&\sum_{k=1}^{\infty}P_{N_{j}}(k)\mathcal{L}_{X_{1}}^{k} (s)\\
=&\mathcal{L}_{X_{1}} (s)\mathcal{L}_{N_{j}}(1-\mathcal{L}_{X_{1}} (s)).
\end{align*}
The stated results now follow from the assumptions.
\end{proof}
\begin{theorem}
$(a)$\,Suppose $N_{1}$ and $N_{2}$  are nabla random variables, which are independent of $X_{i}'$s.  $N_{1}\leq_{D^{^{i}}-Lt-r}N_{2}$ if and only if $X_{i:N_{2}}\leq_{lr}X_{i:N_{1}},\, i=1,2,...,\rho(k).$\\
 $(b)$\,Suppose $N_{1}$ and $N_{2}$  are delta random variables, which are independent of $X_{i}'$s. Then, we have
  $ N_{1}\leq_{D^{^{\gamma}}-Lt-r}N_{2}$ if and only if $X_{\gamma:N_{2}}\leq_{lr}X_{\gamma:N_{1}}.$
\end{theorem}
\begin{remark}
In the recent theorem,  
$(a) $\,in the special case $i=1, $ obviously $P(N_{1}\geqslant 1)=P(N_{2}\geqslant 1)=1$ and so $N_{1}\leq_{d-Lt-r}N_{2}$ if and only if $X_{1:N_{2}}\leq_{lr}X_{1:N_{1}};$\\
 $(b) $\,in the special case $\gamma=\alpha-1,\, \alpha\geq 1, $ obviously $P(N_{1}\geqslant \alpha-1)=P(N_{2}\geqslant \alpha-1)=1$ and then we have
  $ N_{1}\leq_{d^{^{(\alpha-1)}}-Lt-r}N_{2}$ if and only if $X_{\alpha-1:N_{2}}\leq_{lr}X_{\alpha-1:N_{1}}.$
\end{remark}

\section{Applications}
\subsection{Frailty models}
In demography and survival analysis, in order to study the unobserved difference in the risk of death of individuals, frailty
models were introduced to evaluate the effect of the variation on the observed hazards (Vaupel et al. (1979)). These models turn
out to play an important role in providing insights and ideas to explain practical phenomena. The Laplace transform  has an important role in the study of frailty models. For instance, the derivatives of the Laplace transforms are used to obtain general results for the power variance function family. In frailty model,
the joint survival
function for a cluster of size $n$ with covariate information $\mathbb{X}= (x_{1}^{t},..., x_{n}^{t}) $
is obtained from the joint conditional survival function by integrating out
the frailty with respect to the frailty distribution as 
\begin{equation*}
S_{\mathbb{X},f}(\mathbb {T}_{n})=\mathcal{L}(H_{\mathbb{X},c}\left( \mathbb {T}_{n}) \right),
\end{equation*}
where  $H_{\mathbb{X},c}(\mathbb {T}_{n})$ is the sum of the cumulative hazards of the $n_{i}$ subjects in cluster
$i, $ and
 $\mathbb {T}_{n}=(t_{1},t_{2},...,t_{n}).$  For
further details, we refer the reader to Duchateau and
 Janssen (2008).
From the joint survival function, the joint density function for a cluster of size $n$ with covariate information $\mathbb{X}= (x_{1}^{t},..., x_{n}^{t})$ 
 can be obtained as
\begin{equation*}
f_{\mathbb{X},f}(\mathbb {T}_{n})=\prod_{j=1}^{n}h_{X_{j},c}(t_{j})(-1)^{n}\mathcal{L}^{(n)}\left( H_{\mathbb{X},c}(\mathbb {T}_{n})\right).
\end{equation*}
These expressions for the joint survival and density functions are useful in
deriving the contribution of different clusters to the marginal likelihood. The
likelihood contribution corresponds to the survival up to the time of censoring
for censored subjects and the density at the event time for the subjects that
experience the event (Duchateau and Janssen (2008)). Therefore, the conditional likelihood contribution of
cluster $i$ is given by
\begin{equation*}
\prod _{j=1}^{n_{i}}h_{x_{ij},c}^{\delta_{ij}}(y_{ij})(-1)^{d_{i}}\mathcal{L}^{(d_{i})}\left( \sum_{j=1}^{n_{i}}H_{_{ij},c}(y_{ij})\right),
\end{equation*}
where   $h_{x_{ij},c}$ is the conditional hazard function for observation $ j,$ $d_{i}=\sum_{j=1}^{n_{i}}\delta_{ij} $ is the number of events in cluster $i,$ and $y_{ij}$ is the minimum of the event time and censoring time.
\begin{theorem}
Let $X_{1} $ and $X_{2}$ be two population random variables, and
let $U_{1}$ and $U_{2}$ be frailty  random variables.   Then,
\begin{align*}
&U_{1}\leq_{d^{^{(n)}}-Lt-r} U_{2}\Leftrightarrow X_{2} \leq_{hr}X_{1},\\
&U_{1}\leq_{d^{^{(n)}}-Lt-r} U_{2}\Leftrightarrow X_{2}\leq_{lr}X_{1}.
\end{align*}
\end{theorem}

\subsection{Reliability}
The hazard rate ordering is stronger than the usual
stochastic order for random variables, which compares
lifetimes with respect to their hazard rate functions.  This ordering  is useful in reliability theory and
survival analysis  owing to the importance of the hazard rate function in these areas. The following theorem  presents a property of the hazard rate ordering of random variable $N_{s}(X).$
\begin{theorem}
Let $X $ and $Y $ be two delta (nabla) random variables, and
let $N_{s}(X )$ and $N_{s}(Y)$ be as described before, associated with nabla (delta) random variables. Then,
\begin{equation*}
\bar{F} _{Y}\leq_{d^{^{(n-1)}}-Lt-r} \bar{F} _{X}\Leftrightarrow N_{s}(Y )\leq_{hr}N_{s}(X).
\end{equation*}
\end{theorem}
\begin{proof}
Using the forms of $R^{X}_{s}(n)$ and $\mathcal{L}^{**}_{X}(s),$ it is easy to verify the result.
\end{proof}
\subsection{Random extremes}
For a sequence of independent and identical  random variables $X_{1},\,X_{2},...$ and a discrete nabla random variable $N,$ which is independent of $X_{i},$ random extremes $X_{1:N}\equiv min_{_{1\leq i\leq N}X_{i}}$ and $X_{N:N}\equiv max_{_{1\leq i\leq N}X_{i}}$ are of potential applications in reliability, transportation, economics,  etc. The next theorem presents a property of some stochastic orders of random extremes.
\begin{theorem}\label{dast} 
Let $X_{1},\,X_{2},\,...$ be a sequence of non-negative independent and identically  distributed random variables. Let $N_{1}$ and $N_{2}$ be two nabla  random variables, which are independent of the $X_{i}.$\\ 

  $(a)\,\,N_{1}\leqslant _{d-Lt-r}N_{2}$ if and only if $X_{1:N_{2}}\leqslant _{lr} X_{1:N_{1}}\,\,\left( X_{N_{1}:N_{1}}\leqslant _{lr} X_{N_{2}:N_{2}}\right);$\\
 
 $(b)\,\,N_{1}\leqslant _{d-Lt-r}N_{2}$ if and only if $\sum_{i=1}^{\rho(N_{1})}X_{i}\leqslant _{d-Lt-r}\sum_{i=1}^{\rho(N_{2})}X_{i};$\\
 
 $(c)$ If $N_{1}\leqslant _{Lt-r}N_{2}$ then $X_{1:N_{1}}\geqslant _{hr}X_{1:N_{2}}\,\,\left( X_{N_{1}:N_{1}}\leq _{hr}X_{N_{2}:N_{2}}\right);$\\

$(d)$ If $N_{1}\leqslant _{r-Lt-r}N_{2}$ then $X_{1:N_{1}}\geqslant _{rh}X_{1:N_{2}}\,\,\left( X_{N_{1}:N_{1}}\leq _{rh}X_{N_{2}:N_{2}}\right).$
\end{theorem}
\begin{proof}
(a) By using  Equation  (\ref{c}), for $j=1,2$ we can write  
\begin{align*}
f_{X_{1:N_{j}}}(x)
&=-f(x)\sum_{k=1}^{\infty}\rho(k)(1- F (x) )^{\rho(k)-1}P_{N_{j}}(k)\\
&=\frac{d}{dx}\mathcal{L}_{N_{j}}(F).
\end{align*}
Since $ f_{X_{1:N_{2}}}(x)/f_{X_{1:N_{1}}}(x)=(d/dx)\mathcal{L}_{N_{2}}(F)/(d/dx)\mathcal{L}_{N_{1}}(F)$ has the common monotone as $\mathcal{\acute{L}}_{N_{2}}(x)/\mathcal{\acute{L}}_{N_{1}}(x), $ the desired result follows immediately. Note that Equation  (\ref{c}) is applicable only for $i=1,2,..., \rho(k),$ so for the case $i=k$ we  get
\begin{align*}
f_{X_{N_{j}:N_{j}}}(x)
& =\frac{d}{dx} \left\lbrace F(x)\sum_{k=1}^{\infty} (1-\bar{F}(x))^{\rho(k) }P_{N_{j}}(k) \right\rbrace \\
& =\frac{d}{dx}\left\lbrace F(x) \mathcal{L}_{N_{j}}( \bar{F} )\right\rbrace
\end{align*}
and the proof runs a similar manner and hence was omitted. In order to prove 
(b),  we can write 
\begin{align*}
\frac{d}{ds}\mathcal{L}_{X_{1}+...+X_{\rho(N_{j})}} (s)=&\sum_{k=1}^{\infty}\rho(k)P_{N_{j}}(k)\mathcal{L}_{X_{1}}^{\rho(k)-1} (s)\mathcal{\acute{L}}_{X_{1}}(s)\\
&=\mathcal{\acute{L}}_{X_{1}}(s)\sum_{k=1}^{\infty}P_{N_{j}}(k)(1-(1-\mathcal{L}_{X_{1}}(s)))^{\rho(k)-1}\\
&=\frac{d}{ds}\mathcal{L}_{N_{j}}(1-\mathcal{L}_{X_{1}}(s)),
\end{align*}
then, for all $s\geq 0,$
\begin{align*}
\dfrac{\frac{d}{ds}\mathcal{L}_{X_{1}+...+X_{\rho(N_{2})}} (s)}{\frac{d}{ds}\mathcal{L}_{X_{1}+...+X_{\rho(N_{1})}} (s)}=\dfrac{\frac{d}{ds}\mathcal{L}_{N_{2}} (1-\mathcal{L}_{X_{1}}(s))}{\frac{d}{ds}\mathcal{L}_{N_{1}} (1-\mathcal{L}_{X_{1}}(s))}.
\end{align*}
To prove (c), note that   
\begin{align*}
\bar{F}_{X_{1:N_{j}}}=\sum_{k=1}^{\infty}(1-F(x))^{\rho(k)} P_{N_{j}}(k)=\mathcal{L}_{N_{j}}(F (x)),\,\,\,\,\,j=1,2.
\end{align*}
In order to prove $X_{1:N_{1}}\geqslant_{hr}X_{1:N_{2}},$ we need to show that $\bar{F}_{X_{1:N_{2}}}(x)/ \bar{F}_{X_{1:N_{1}}}(x)$ 
is decreasing in $x\geqslant 0,$ but this follows from the fact that
\begin{align*}
\frac{\bar{F}_{X_{1:N_{2}}}(x)}{\bar{F}_{X_{1:N_{1}}}(x)}=\frac{\mathcal{L}_{N_{2}}(F(x))}{\mathcal{L}_{N_{1}}(F(x))},
\end{align*}
and from $N_{1}\leqslant _{Lt-r}N_{2}.$  Similarly, we have $\bar{F}_{X_{N_{j}:N_{j}}}=\bar{F}(x)\mathcal{L}_{N_{j}}(F ).$
The proof of (d) is similar. 
\end{proof}

\section*{Acknowledgment}
The authors are grateful to the referees for the comments and suggestions which led
to substantial improvement of the manuscript. 
The work was mainly done when Fatemeh Gharari was a visiting student in McMaster University, Hamilton, Canada. She thanks Professor Narayanaswamy Balakrishnan for assistance in preparation of the manuscript.


\end{document}